\documentstyle [12pt,twoside]{article}
\font\twelvebb=msym10   at 12pt
\font\twelvegoth=eufm10    at 12pt
\newfam\gothfam
\textfont\gothfam=\twelvegoth
\def\goth{\fam\gothfam\twelvegoth}
\newfam\bbfam
\textfont\bbfam=\twelvebb
\def\bb{\fam\bbfam\twelvebb}

\newcommand{\nc}{\newcommand}
\nc{\la}{\lambda}
\nc{\ga}{\gamma}
\nc{\Ga}{\Gamma}
\nc{\al}{\alpha}
\nc{\be}{\beta}
\nc{\de}{\delta}
\nc{\De}{\Delta}
\nc{\nf}{\infty}
\nc{\ra}{\rightarrow}
\nc{\Ra}{\Rightarrow}
\nc{\Lra}{\Longrightarrow}
\nc{\beq}{\begin{equation}}
\nc{\eeq}{\end{equation}}
\nc{\beqa}{\begin{eqnarray}}
\nc{\eeqa}{\end{eqnarray}}
\nc{\bfi}{\begin{figure}}
\nc{\efi}{\end{figure}}
\nc{\prh}[2]{\left(\hspace{-0.2cm}\begin{array}{c}#1
\\#2 \end{array}\hspace{-0.2cm}\right)}
\nc{\prf}[3]{\left(\hspace{-0.2cm}\begin{array}{c}#1
\\#2 \end{array}\hspace{-0.1cm};#3\right)}
\nc{\cro}[2]{\left[\begin{array}{c}#1 \\#2
\end{array}\right]}
\nc{\fx}[5]{{_{#1}F_{#2}}\prf{#3}{#4}{#5}}
\nc{\ff}[4]{{_{#1}F_{#2}}\prh{#3}{#4}}

\title{\bf Some results on co-recursive associated Laguerre and Jacobi
polynomials}

\author{\bf Jean  Letessier\thanks{\noindent Laboratoire de Physique
Th\'eorique et
Hautes Energies, Unit\'e
associ\'ee au CNRS UA 280,~Universit\'e PARIS 7, Tour 24,
5\`e \'et., 2 Place Jussieu, F-75251 CEDEX 05.}}

\date{November 1992}   

\begin{document}           

\maketitle                 
\centerline{ To be published in the special issue of the SIAM J. Math. Anal.}
\centerline{ dedicated to R. Askey and F. Olver.}

\pagenumbering {arabic}
\begin{abstract}   We present  results on co-recursive  associated Laguerre and
Jacobi
polynomials which are of interest for the solution of the Chapman-Kolmogorov
equations of some birth and death processes with or without absorption.
Explicit forms,
generating functions, and absolutely continuous part of the spectral measures
are given.
We derive fourth-order differential equations satisfied by the polynomials
 with a special attention to some simple limiting cases.
\vspace{2mm}

\vspace{2mm}

\noindent{\bf AMS(MOS) subject classifications.} 33A65, 60J80, 33A30.
\end{abstract}
\section {Introduction}

Starting from a sequence of orthogonal polynomials $\{P_n\}_{n\geq0}$ defined
by the recurrence relation
\beq\label{RR1}
P_{n+2}(x)=(x-\be_{n+1})P_{n+1}(x)-\ga_{n+1}P_n(x),\ n\geq 0,
\eeq
and the initial conditions
\beq\label{IN}
P_0(x)=1,\hspace{1cm}P_1(x)=x-\be_0,
\eeq
with $\be_n$, $\ga_n\in{\bb C}$  and $\ga_n\neq 0$,
several modifications were  considered:
\begin{itemize}
\item Associated polynomials arise when we replace $n$ by $n+c$ in the
coefficients $\be_n$ and $\ga_n$ (keeping $\ga_n\neq 0$).
If $c$ is an integer
$k$ these polynomials are called associated of order $k$.
The associated
polynomials of order one are the numerator polynomials.
\item Co-recursive
polynomials arise when we replace $\be_0$ by $\be_0+\mu$.
 \item Perturbed
polynomials arise when we replace $\ga_1$ by $\la\ga_1$, ($\la>0$).
 \item Co-recursive of this perturbed polynomials arise when the two previous
modifications are made together.
\item Generalized co-recursive and perturbed
polynomials arise when we change $\be_n$ and/or $\ga_n$ at any level $n$.
\end{itemize}

In the study of birth and death processes orthogonal polynomials,
in particular all the hypergeometric
families of the Askey scheme \cite{AND85,LAB85} and their corresponding
associated families, play a
primordial role in the Karlin-McGregor solution of the Chapman-Kolmogorov
equation \cite{KAR58,ISM90C}
\beq
p_{mn}(t)\sim \int_0^\nf e^{-xt}P_m(x)P_n(x)d\phi(x).
\eeq
In certain birth and death processes, zero-related polynomials
\cite{ISM88,ISM89,ISM90}
arise in a natural way. Zero-related polynomials are special  co-recursive
associated
polynomials.

More generally co-recursive and generalized co-recursive polynomials are
involved in the solution of the Chapman-Kolmogorov equation of birth and death
processes
with absorption or killing \cite{ISM90C,KAR82,KAR82B}.

The purpose of this paper is to present some results on co-recursive associated
Laguerre and Jacobi polynomials which are of special interest in the study of
birth and death processes with linear and rational rates respectively.
The Laguerre polynomial families are involved in processes for which the birth
and death
rates are of the form \cite{ISM88}
\beq
\la_n=n+\al+c+1,\hspace{0.5cm}\mu_{n+1}=n+c+1,\hspace{0.5cm}n\ge 0,
\hspace{1.5cm}\mu_{0}=c-\mu.
\eeq
The Jacobi polynomial families are involved when the rates are of the form
\beqa
\la_n&&\hspace{-0.6cm}=\frac{2(n+c+\al+\be+1)(n+c+\be+1)}
{(2n+2c+\al+\be+1)(2n+2c+\al+\be+2)},\hspace{0.5cm}n\ge0,\\
\mu_n&&\hspace{-0.6cm}
=\frac{2(n+c)(n+c+\al)}{(2n+2c+\al+\be)(2n+2c+\al+\be+1)},
\hspace{1.2cm}n>0,\\
\mu_0&&\hspace{-0.6cm}=\frac{2c(c+\al)}{(2c+\al+\be)(2c+\al+\be+1)}-\mu.
\eeqa

In both cases $\mu_0=0$ correspond to the ``honest''  \cite{REU57}
linear processes (i.e. processes for which the sum of the probabilities
$p_{mn}(t)$
is equal to 1).
Cases $\mu_0=$ Const. $\neq0$ correspond to processes with absorption
and are not ``honest''.
However if $\mu_0=c$ in the Laguerre case or
$\mu_0=\frac{2c(c+\al)}{(2c+\al+\be)(2c+\al+\be+1)}$ in
the Jacobi case, the corresponding processes are simply solved using associated
polynomials.

In  section \ref{SLC} we explain the method used by applying it
to the Laguerre case.
In \ref{SSLEF} we give an explicit expression for the
co-recursive associated Laguerre (CAL) polynomials. In \ref{SSCALM} we derive
a generating function of them and we found the absolutely continuous part of
the
spectral measure.
The subsection \ref{SSLDE} is devoted to the derivation, using the Orr's
method,  of a
fourth-order differential equation satisfied by the CAL polynomials.
In \ref{SPC} we present results obtain in some limiting cases
among which a new simple case of associated Laguerre polynomials.

In section \ref{SJC} we present briefly some results corresponding to the
Jacobi case.
In \ref{SSJEF} we give an explicit expression for the co-recursive associated
Jacobi (CAJ)
polynomials.
In \ref{SSCAJGF} we present a generating function and in \ref{SSCAJM} we give
the
absolutely continuous component of the spectral measure.
The subsection \ref{SSJPC} is devoted to some limiting cases of CAJ
polynomials
for which we give fourth-order differential equations they satisfy.
In section \ref{SCONC} are some concluding remarks.

We use the notation of \cite{HTF1} for the special functions used in this work.
We don't give validity conditions on
the parameters of the used hypergeometric functions, analytic continuations or
limiting processes giving, in general, valid formulas.
We use Slater's notation \cite{SLA66} for the product of $\Ga$ functions
\beq
\Ga\prh{\al_1,\ldots,\al_p}{\be_1,\ldots,\be_q}
=\left.{\displaystyle\prod_{i=1}^p\Ga(\al_i)}\right/{\displaystyle
\prod_{i=1}^q\Ga(\be_i)}.
\eeq

\section{The case of Laguerre polynomials}\label{SLC}

Replacing $n$ by $n+c$ in the recurrence relation of the Laguerre polynomials
we obtain the recurrence relation satisfied by  the associated Laguerre
polynomials $L^\al_{n}(x;c)$
\beq\label{RRL}
(2n+2c+\al+1-x)p_n=(n+c+1)p_{n+1} +(n+c+\al)p_{n-1}.
\eeq
To complete the definition of the polynomials $L^\al_{n}(x;c)$ the initial
condition
\beq\label{ICL}
L^\al_{-1}(x;c)=0,\hspace{1cm}L^\al_{0}(x;c)=1,
\eeq
has to be imposed. These polynomials are orthogonal with respect to a positive
measure when $(n+c)(n+\al+c)>0,\ \forall n>0$. (See \cite{ASK84B} for details).

Note that if we consider the monic polynomials ${\bf L}^\al_{n}(x;c)$ defined
by
\beq
{\bf L}^\al_{n}(x;c)=(-1)^n(c+1)_nL^\al_{n}(x;c)
\eeq
they satisfy the recurrence relation
\beq\label{RRML}
(x-2n-2c-\al-1){\bf L}^\al_{n}(x;c)={\bf L}^\al_{n+1}(x;c)+(n+c)(n+c+\al){\bf
L}^\al_{n-1}(x;c).
\eeq
We can see that this recurrence is invariant in the
transformation ${\cal T}$ defined by
\beq\label{TT}
{\cal T}(c,\al)=(c+\al,-\al).
\eeq

\subsection{Explicit representation of the CAL polynomials}  \label{SSLEF}
Equations (\ref{RRL}) and (\ref{ICL}) give for $L^\al_{1}(x;c)$
\beq
L^\al_{1}(x;c)=-\frac{1}{c+1}(x-2c-\al-1)
\eeq
The CAL polynomials $L^\al_{n}(x;c,\mu)$
satisfy the same recurrence relation (\ref{RRL}) with a shift $\mu$ on the
monic polynomial of first degree, i.e.
\beq \label{L1}
L^\al_{1}(x;c,\mu)=-\frac{1}{c+1}(x+\mu-2c-\al-1).
\eeq
To obtain $L^\al_{1}(x;c,\mu)$ with (\ref{RRL}) we have to impose the initial
condition for the CAL polynomials
\beq\label{ICAL}
L^\al_{-1}(x;c,\mu)=\frac{\mu}{c+\al},\hspace{1cm}L^\al_{0}(x;c,\mu)=0.
\eeq
Even for $c+\al\ra 0$ this initial condition in the recurrence relations
(\ref{RRL}) leads to
(\ref{L1}).

We know two linearly independent solutions of (\ref{RRL}) in terms of confluent
hypergeometric functions:
\beq
\hspace{1cm}u_n=\frac{(c+\al+1)_n}{(c+1)_n}\fx{1}{1}{-n-c}{1+\al}{x}
\hspace{1cm}{\rm and}\hspace{1cm}
v_n=\fx{1}{1}{-n-c-\al}{1-\al}{x}.
\eeq
Writing the polynomials $L^\al_{n}(x;c,\mu)$ as a linear combination
\beq\label{LC}
L^\al_{n}(x;c,\mu)=Au_n+Bv_n
\eeq
and using the initial condition (\ref{ICAL}) we obtain
\beq\label{AA}
A=\frac{1}{\Delta}\left[\frac{\mu}{c+\al}v_0-v_{-1}\right]\hspace{1cm}{\rm
and}\hspace{1cm}
B=-\frac{1}{\Delta}\left[\frac{\mu}{c+\al}u_0-u_{-1}\right],
\eeq
where $\Delta$ may be calculated using contiguous relations of confluent
hypergeometric functions \cite[page 253--254]{HTF1}
\beq
\Delta=u_{-1}v_0-u_0v_{-1}=-\frac{\al}{c+\al}e^x.
\eeq
With the help of the relation
\beq\label{2F2}
\ga\fx{1}{1}{1-\ga}{\al}{x}-\be\fx{1}{1}{-\ga}{\al}{x}=(\ga-\be)
\fx{2}{2}{-\ga,\be-\ga+1}{\al,\be-\ga}{x},
\eeq
we can write the CAL polynomials as
\beqa
L^\al_{n}(x;c,\mu)=\frac{e^{-x}}{(c+1)_n}&&\hspace{-0.6cm}(1+{\cal
T})\nonumber\\
&&\hspace{-0.6cm}\hspace{-0.9cm}\times\frac{\mu-c}{\al}
{(c+1)_n}\fx{2}{2}{-c,\mu-c+1}{1+\al,\mu-c}{x}\fx{1}{1}{-n-c-\al}{1-\al}{x}.
\label{CALP}\eeqa
This representation is the one we use to derive the fourth-order differential
equation in section \ref{SSLDE}.
It is valid only for $\al\neq 0,\pm 1,\pm 2,\ldots$ but these
restrictions can be removed by limiting processes.

Following the same way as in \cite{ASK84B} we can find an explicit
representation.
 We first transform the ${_1F_1}$ in (\ref{AA}) using Kummer's transformation
\cite[page 253]{HTF1}
\beq\label{K1F1}
\fx{1}{1}{a}{c}{x}=e^x\fx{1}{1}{c-a}{c}{-x}.
\eeq
 Then we use formula
 \beq
 \fx{1}{1}{a}{b}{x}\fx{1}{1}{c}{d}{-x}=\sum_{k=0}^\nf
 \frac{(a)_kx^k}{k!(b)_k}\fx{3}{2}{-k,1-k-b,c}{1-k-a,d}{1},
 \eeq
 for the four products in (\ref{LC}).

 Applying the three term relation \cite[Eq. (2) page 15]{BAI72} to the
four resulting ${_3F_2}(1)$ gives a sum of two terms which we can group to
obtain
the explicit form
 \beqa
L^\al_{n}(x;c,\mu)=\frac{(c+\al+1)_n}{n!}\sum_{k=0}^n&&\hspace{-0.6cm}
\frac{(-n)_k}{(c+1)_k(c+\al+1)_k}\nonumber\\
&&\hspace{-0.6cm}\hspace{0.6cm}\times\fx{4}{3}{k-n,c,c+\al,G+1}
{c+k+1,c+\al+k+1,G}{1}x^k,\label{EFLP}
\eeqa
where
\beq
G=\frac{c(c+\al)}{\mu}.
\eeq
Take care of the limiting processes when $n=k$ and when $c$ or $c+\al=0$.
The representation (\ref{EFLP}) can also be proved using
the generating function (\ref{GFL}).

\subsection{Generating function, spectral measure}    \label{SSCALM}

Let $F(x,w)$ be a generating function of the CAL polynomials
$L^\al_{n}(x;c,\mu)$
\beq
F(x,w)=\sum_{n=0}^\nf w^nL^\al_{n}(x;c,\mu).
\eeq
The recurrence relation (\ref{RRL}) and the initial condition (\ref{ICAL})
lead to the following differential equation for $F(x,w)$
\beq\label{DEGF}
\hspace{1cm}
w(1-w)^2\frac{\partial\ }{\partial
w}F(x,w)+[(1-w)(c-(c+\al+1)w)+xw]F(x,w)=c-\mu
w. \eeq
The function $F(x,w)$ is normalized by the condition $F(x,0)=1$ and due to the
orthogonality of the $L^\al_{n}(x;c,\mu)$, we have the boundary condition
\beq\label{BC} \int_0^\nf F(x,w)d\phi(x)=1,
\eeq
where $d\phi(x) $ is the spectral measure.

The solution of the differential equation (\ref{DEGF}) which is bounded at
$w=0$ is easily obtain, for $c>0$, following the same methode as in
\cite{ISM88}
\beqa
F(x,w)=w^{-c}(1-w)^{-\al-1}&&\hspace{-0.6cm}
\exp{\left[-\frac{x}{1-w}\right]}\nonumber\\
&&\times\int_0^wu^{c-1}(1-u)^{\al-1}(c-\mu
u)\exp\left[\frac{x}{1-u}\right]du.\label{GFL}
\eeqa
Changing variables according to
\beq
u=\frac{\tau}{1+\tau},\hspace{1cm}w=\frac{z}{1+z},
\eeq
and integrating both side of (\ref{GFL}) with respect to $d\phi(x)$, taking
into account (\ref{BC}),  leads to
\beq
z^c(1+z)^{-1-\al-c}=\int_0^\nf\int_0^z\tau^{c-1}(1+\tau)^{-1-\al-c}
[c+\tau(c-\mu)]\exp{[-x(z-\tau)]d\tau d\phi(x)}.
\eeq

Taking the Laplace transform of the above identity we obtain for the Stieltjes
transform of the measure $d\phi(x)$ the relationship
\beq\label{STL}
s(p)=\int_0^\nf\frac{d\phi(x)}{x+p}=\frac{\Psi(c+1,1-\al;p)}
{\Psi(c,-\al;p)+(c-\mu)\Psi(c+1,1-\al;p)},
\eeq
which we can rewrite formaly, using the  expression of the Tricomi function
$\Psi$ in
terms of generalized hypergeometric functions \cite[page 257]{HTF1}, on the
form
\beq
s(p)=p^c{\Psi(c+1,1-\al;p)}{\fx{3}{1}{c,c+\al,\frac{c(c+\al)}{\mu}+1}
{\frac{c(c+\al)}{\mu}}{-\frac{1}{p}}}^{-1},
\eeq
where the principal branch of the ${_3F_1}$ is considered as that one of the
$\Psi$
function. The function $\Psi (a,b;p)$ having no zeros for $|\arg p|\leq\pi$
the denominator of (\ref{STL}) has no zeros in this region at least for
$\mu\leq c$.

The CAL polynomials belong to the Laguerre-Hahn family of orthogonal
polynomials
and are of class zero \cite{MAR91}. It is
easily verified that the Stieltjes transform $s(p)$ of the measure,
 calculated in (\ref{STL}),
is a solution of the Riccati equation
\beq
ps'(p)=\left[\mu\,p-(c-\mu)(\al+c-\mu)\right]s^2(p)
+\left[p+\al+2(c-\mu)\right]s(p)-1.
\eeq

The absolutely continuous part of the measure $d\phi(x)$ can be computed using
the Perron-Stieltjes inversion formula. Details of the method are in
\cite{ISM79} and we obtain
\beq
\phi'(x)=\frac{1}{\Ga(c+1)\Ga(c+\al+1)}\frac{x^{\al}e^{-x}}
{\left|\Psi(c,-\al;xe^{i\pi})+(c-\mu)\Psi(c+1,1-\al;xe^{i\pi})\right|^2},
\eeq
we can write formaly
\beq
\phi'(x)=\frac{x^{\al+2c}e^{-x}}{\Ga(c+1,c+\al+1)}
\left|\fx{3}{1}{c,c+\al,\frac{c(c+\al)}{\mu}+1}
{\frac{c(c+\al)}{\mu}}{-\frac{e^{i\pi}}{x}}\right|^{-2}.
\eeq

The CAL polynomials $L_n^\al(x;c,\mu)$ satisfy the orthogonality relation,
valid
at least for $\mu\leq c$, $c\geq0$, $\al+c>-1$,

\beq
\int_0^\nf L_n^\al(x;c,\mu)L_m^\al(x;c,\mu)d\phi(x)=
\frac{(c+\al+1)_n}{(c+1)_n}\delta_{mn}.
\eeq

\subsection{Fourth-order differential equation}        \label{SSLDE}

The CAL polynomials  verify a fourth-order differential equation
\cite{ATK81,MAR91}.
One way to obtain this fourth-order differential equation
 is to start from their explicit form (\ref{CALP}). The righthand
side of (\ref{CALP}) is a sum of two products of a ${_1F_1}$ times a ${_2F_2}$.
The ${_1F_1}$ are solution of a second-order differential equation but
 for the ${_2F_2}$ a third-order one is expected.
 In fact  the ${_2F_2}$ involved in (\ref{CALP}) are of the form
\beq
y(x)=\fx{2}{2}{b,e+1}{d,e}{x},
\eeq
 which can be shown with little effort to be  a solution of the second-order
differential
 equation
\beqa
x[(b-e)x&&\hspace{-0.6cm}+e(d-e-1)]y''(x)-\{(b-e)x^2+[e(2d-e-2)+b(1-d)]x\\
&&\hspace{-0.6cm}\hspace{0.4cm}-ed(d-e-1)\}y'(x)
-b[(b-e)x+(e+1)(d-e-1)]y(x)=0.\nonumber
\eeqa

So we can use the Orr method to obtain the differential equation satisfied by
the products involved in (\ref{CALP}) \cite{ORR00}.
Let us notice that the second product being obtained by the transformation
$\cal T$ of the first one, the fourth-order differential equation will have to
be invariant under this transformation.

The function $y=e^{-x}\fx{1}{1}{-n-c-\al}{1-\al}{x}$ is solution of
\beq\label{DIF1}
xy''+(1-\al+x)y'+(1+n+c)y=0
\eeq
and the function $z=\fx{2}{2}{-c,\mu-c+1}{1+\al,\mu-c}{x}$ of
\beqa
x[\mu x&&\hspace{-0.6cm}+(c-\mu)(c+\al-\mu)]z''(x)
-\{\mu x^2+[(c-\mu)(c+\al-\mu)-\al\mu]x\\
&&\hspace{-0.6cm}\hspace{-0.8cm}-(\al+1)(c-\mu)(c+\al-\mu)\}z'(x)
+c[(\mu x+(c-\mu+1)(c+\al-\mu)]z(x)=0.\nonumber
\label{DIF2}
\eeqa
Changing the functions $y$ and $z$ to $y=fv$ and $z=gw$
with
\beqa
f&&\hspace{-0.6cm}=x^{\frac{\al-1}{2}}e^{\frac{x}{2}},\\
g&&\hspace{-0.6cm}=[(c-\mu)(c+\al-\mu)+\mu x]^\frac{1}{2}
x^{-\frac{\al+1}{2}}e^{\frac{x}{2}},
\eeqa
we obtain the normal form of the differential equations (\ref{DIF1})
and (\ref{DIF2})
\beq
v''+Iv=0,\hspace{1cm}w''+Jw=0.
\eeq
The product $u=vw$ is solution of the fourth-order differential equation (see
\cite[page 146]{WAT44})
\beq
\frac{d\ }{dx}\left[\frac{u'''+2(I+J)u'+(I'+J')u}{I-J}\right]=-(I-J)u
\eeq
Finally we obtain the needed equation for the CAL polynomials setting
$y(x)=fgu$.

Details of this calculations are very difficult to write explicitly and were
achieved with the help of the MAPLE computer algebra \cite{MAP88}. Although
with
this help the fourth-order differential equation for the CAL polynomials
is not easy to find. We give it as a {\em curiosity\/}:
\beq\label{CALDE}
c_4 y^{(4)}(x)+c_3 y^{(3)}(x)+c_2 y^{(2)}(x)+c_1 y^{(1)}(x)+c_0 y(x)=0,
\eeq
with\vspace{-0.3cm}
\beqa
c_4=&&\hspace{-0.6cm}x^2(2 A  x^2+B x+2 C),\\
c_3=&&\hspace{-0.6cm}2 x(3 A  x^2+2 B x+5 C),\\
c_2=&&\hspace{-0.6cm}-2 A  x^4+D  x^3+E  x^2+F x+G,\\
c_1=&&\hspace{-0.6cm}-4 A  x^3+H  x^2-4 A  x^3+I x+J,\\
c_0=&&\hspace{-0.6cm}n(n+1)(2 A  x^2+K x+2 L),
\eeqa
where
\beqa
A=&&\hspace{-0.6cm}\mu^2(1+2 n),\nonumber\\
B=&&\hspace{-0.6cm}\mu(2(1+4 n)\mu^2-(4(1+2 n)(2 c+\al)-1)\mu+2 c(3+4
n)(c+\al)),\nonumber\\
C=&&\hspace{-0.6cm}(c-\mu)(c+\al-\mu)(2 n\mu^2-((1+2 n)(2
c+\al)+1)\mu+2 c(n+1)(c+\al)),\nonumber\\
D=&&\hspace{-0.6cm}-\mu(2(1+4 n)\mu^2-(2(1+2
n)(8 c+4\al+1+2 n)-1)\mu+2 c(3+4 n)(c+\al)), \nonumber\\
E=&&\hspace{-0.6cm}-4n\mu^4\!  +\! ((1+4 n)(4n+12 c+6\al+3)+3)\mu^3
\! -\! 2(1+2 n)((2 c+\al)\nonumber\\
&&\hspace{-0.6cm}\times(4n+12 c+6\al+2)-2c(c+\al)-1)\mu^2 +c(c+\al)((3+4
n)\nonumber\\
&&\hspace{-0.6cm}\times(4 n+12 c+6\al+1)+3)\mu-4 c^2(c+\al)^2(n+1),\nonumber\\
F=&&\hspace{-0.6cm}8 n(2 c+\al+n+1)\mu^4\nonumber\\
&&\hspace{-0.6cm}-((1+4 n)((2
c+\al)(4 n+12 c+6\al+3) -8  c(c+\al)-1/2)+10 c+5\al+5/2)\mu^3 \nonumber\\
&&\hspace{-0.6cm}+((1+4 n)((\al^2+6 c(c+\al)(2
n+8c+4\al+3/2)-(2c+\al)(12c(c+\al)
+1/2)) \nonumber\\
&&\hspace{-0.6cm}+(2\al^2+6 c(c+\al)-1/4)(4 c+2\al+23/6)
-25/6\al^2 +71/24)\mu^2\nonumber\\
&&\hspace{-0.6cm}-2 c(c+\al)((1+4 n)((2 c+\al)(2 n+4 c+2\al+3/2)
+2 \al^2-1/4)\nonumber\\
&&\hspace{-0.6cm}+(2 c+\al)(8 c+4\al+5/2)+2\al^2-7/4)\mu+8
c^2(c+\al)^2(n+1)(n+2
c+\al),\nonumber\\
G=&&\hspace{-0.6cm}-2(\al-2)(\al+2)(c-\mu)(c+\al-\mu)(2(c-\mu)
(c+\al-\mu)n\nonumber\\
&&\hspace{-0.6cm}-(2c+\al+1)\mu+2 c(c+\al)),\nonumber\\
H=&&\hspace{-0.6cm}2 \mu(-2(5 n+1)\mu^2
+((1+2 n)(n+12c+6\al +1/2)-3/2)\mu-2c(c+\al)(5 n+4)),\nonumber\\
 I=&&\hspace{-0.6cm}-12n\mu^4 +2((n+1/5)(8n+40c+20\al+7/5)+93/25)\mu^3\nonumber
\\  &&\hspace{-0.6cm}+2(-(1+2 n)((2 c+\al)(4
n+12 c+6\al+2)+10c(c+\al)-1)-6 c-3\al)\mu^2 \nonumber \\
&&\hspace{-0.6cm}+2c(c+\al)((4/5+n)(
 8 n+40 c+20\al+33/5)+93/25)\mu-12  c^2(c+\al)^2(n+1),\nonumber\\
J=&&\hspace{-0.6cm}4(\mu-c)(\mu-c-\al)(3 n(n+2 c+\al+1)\mu^2+(-1/4(1+2 n)
((2 c+\al)(6 n+8c \nonumber\\
&&\hspace{-0.6cm}+4\al+3)+8c(c+\al)+8)-9/2 c-9/4\al)\mu+3 c(c+\al)(n+1)(n+2
c+\al)),\nonumber\\
K=&&\hspace{-0.6cm}\mu(2(4 n-1)\mu^2-(4(1+2 n)(2 c+\al)-3)\mu+2 c(4
n+5)(c+\al)),\nonumber\\
 L=&&\hspace{-0.6cm}(c-\mu)(c+\al-\mu)(2(n-1)\mu^2-((1+2 n)(2
c+\al)+6)\mu+2 c(n+2)(c+\al)). \nonumber
\eeqa
 Note the invariance of the differential equation
(\ref{CALDE}) by the transformation ${\cal T}$ defined in (\ref{TT}).

\subsection{Particular cases}\label{SPC}
We now give the different results corresponding to limiting cases of special
interest.

\subsubsection{Limit $c=0$}   \label{SSCLP}
In this limit we obtain from (\ref{CALP}) the co-recursive Laguerre
polynomials.
\beqa
L_n^\al(x;0,\mu)=&&\hspace{-0.6cm}
\frac{e^{-x}}{\al}\left\{\frac{(\al-\mu)(\al+1)_n}{n!}
\fx{2}{2}{-\al,\mu-\al+1}{1-\al,\mu-\al}{x}\fx{1}{1}{-n}{1+\al}{x}
\right.\nonumber\\
&&\hspace{-0.6cm}\hspace{7.4cm}\left.+\mu\fx{1}{1}{-n-\al}{1-\al}{x}\right\}.
\label{CL}\eeqa
The limit $\mu=0$ in (\ref{CL}) leads back to the classical Laguerre
polynomials.  An explicit form is
\beqa
L^\al_{n}(x;0,\mu)=\frac{(\al+1)_n}{n!}&&\hspace{-0.6cm}\sum_{k=0}^n
\frac{(-n)_k}{k!(1+\al)_k}\nonumber\\
&&\hspace{-0.6cm}\times\left[1+\frac{\mu(k-n)}{(1+k)(1+\al+k)}
\fx{3}{2}{1+k-n,1+\al,1}{k+\al+2,k+2}{1}\right]x^k,
\eeqa
and the corresponding absolutely continuous part of the measure is given
for $\mu\leq 0$, $\al>-1$, by
\beq
\phi'(x)=\frac{x^\al
e^{-x}}{\Ga(1+\al)}\left|1-\mu\Psi(1,1-\al;xe^{-i\pi})\right|^{-2},
 \eeq
where the limit $\mu=0$ is also straightforward.

It is easy to see that the differential equation (\ref{CALDE}) satisfied by the
co-recursive Laguerre polynomials can be factorized in the limit $c=0$ to
obtain the fourth-order factorized (2+2) differential equation
\beqa
&&\hspace{-0.6cm}\hspace{-0cm}\left[xA(x){\rm D}^2+\{(2+\al-x)A(x)-xB(x)\}{\rm
D}+(n-1)A(x)+(x-\al-1)B(x)+C(x)\right]\nonumber\\
&&\hspace{-0.6cm}\hspace{5.6cm}\times\left[x  {\rm D}^2+(x+1-\al){\rm
D}+n+1\right]L_n^\al(x;0,\mu)=0,
\eeqa
where ${\rm D}\,\equiv\,d/dx$ and
\beqa
A(x)&&\hspace{-0.6cm}=3 x+2(x-\al+\mu)\{2 n(x-\al+\mu)+x-\al-1\},\\
B(x)&&\hspace{-0.6cm}=1+2 x-2 \al+2(1+4 n)(x-\al+\mu),\\
C(x)&&\hspace{-0.6cm}=(1+2 x-2 \al)\{1+\al-x-2 n(x-\al+\mu)\}+3 x(1+4 n).
\eeqa
The comparison with the differential equation given in \cite[Eq. 34--35]{RON89}
requires some
attention because of a few misprints.

\subsubsection{Limit $c=-\al$}

In this limit we obtain a special class of CAL polynomials corresponding to the
associated
Laguerre polynomials for which
\beq\label{ALM}
L_n^\al(x;-\al)=\frac{n!}{(1-\al)_n}L_n^{-\al}(x).
\eeq
We can write the $L_n^\al(x;-\al,\mu)$
\beqa
L_n^\al(x;-\al,\mu)=&&\hspace{-0.6cm}\frac{n!}{(1-\al)_n}\frac{e^{-x}}{-\al}
\left\{\frac{(-\al-\mu)(1-\al)_n}{n!}
\fx{2}{2}{\al,\mu+\al+1}{1+\al,\mu+\al}{x}\fx{1}{1}{-n}{1-\al}{x}
\right.\nonumber\\
&&\hspace{-0.6cm}\hspace{7.5cm}\left.+\mu\fx{1}{1}{-n+\al}{1+\al}{x}\right\}.
\label{CLM}\eeqa
which gives (\ref{ALM}) in the limit $\mu =0$. Except for the
global factor $\frac{n!}{(1-\al)_n}$,
 (\ref{CLM}) is obtained from (\ref{CL}) changing $\al$ to $-\al$.
The corresponding measure and
differential equation are obtained in the same way from \ref{SSCLP}.  An
explicit form is
\beqa
L^\al_{n}(x;-\al,\mu)=\sum_{k=0}^n&&\hspace{-0.6cm}
\frac{(-n)_k}{k!(1-\al)_k}\nonumber\\
&&\hspace{-0.6cm}\times\left[1+\frac{\mu(k-n)}{(1+k)(1-\al+k)}
\fx{3}{2}{1+k-n,1-\al,1}{k-\al+2,k+2}{1}\right]x^k.
\eeqa

\subsubsection{Limit $\mu=0$}
In this limit we obtain the associated Laguerre polynomials studied in
\cite{ASK84B} and \cite{ISM88}
\beq
L_n^\al(x,c)=\frac{e^{-x}}{(c+1)_n}(1+{\cal T})\frac{(c+\al)_{n+1}}{\al}
\fx{1}{1}{1-c-\al}{1-\al}{x}\fx{1}{1}{-n-c}{1+\al}{x},
\eeq
with the measure
\beq
\phi'(x)=\frac{x^\al e^{-x}}{\Ga(1+c,1+c+\al)}
\left|\Psi(c,1-\al;xe^{-i\pi})\right|^{-2}.
\eeq
An explicit form is
\beqa
L^\al_{n}(x;c)=\frac{(c+\al+1)_n}{n!}\sum_{k=0}^n&&\hspace{-0.6cm}
\frac{(-n)_k}{(c+1)_k(c+\al+1)_k}\nonumber\\
&&\hspace{-0.6cm}\times\fx{3}{2}{k-n,c,c+\al}{c+k+1,c+\al+k+1}{1}x^k,
\eeqa
The limit $c=0$ leads back to the Laguerre polynomial case.
The coefficients of the differential equation (\ref{CALDE}) satisfied by the
associated Laguerre polynomials are now very simple
\beq
\hspace{0.9cm}c_4 =  x^2,\hspace{0.2cm}c_3 =  5 x,\hspace{0.2cm}
c_2 =-x(x-2 F)-\al^2+4,\hspace{0.2cm}
c_1 =   3 (F- x),\hspace{0.2cm}c_0 =  n(n+2),
\eeq
with $F=n+2c+\al$.
This differential equation was first given by Hahn \cite[Eq. 22]{HAH40}.
 See \cite{RON88,RON89}
for the special factorizable case $c=1$ and \cite{BEL91,RON91} when $c$ is an
integer.

\subsubsection{Limit $\mu=c$ }\label{ZRL}
In this limit we obtain the so called {\em zero related\/} Laguerre polynomials
studied in  \cite{ISM88}. Note  the symmetry ${\cal T}$ of the monic
polynomials is now broken.
\beqa
{\cal L}_n^\al(x,c)=e^{-x}&&\hspace{-0.6cm}\left\{\frac{
(c+\al+1)_{n}}{(c+1)_n}
\fx{1}{1}{-c-\al}{-\al}{x}\fx{1}{1}{-n-c}{1+\al}{x}\right.\nonumber\\
&&\hspace{-0.6cm}-\frac{c}{\al(\al+1)}\ x
\left.\fx{1}{1}{1-c}{2+\al}{x}\fx{1}{1}{-n-c-\al}{1-\al}{x}\right\},
\eeqa
and the measure
\beq
\phi'(x)=\frac{x^\al e^{-x}}{\Ga(1+c,1+c+\al)}
\left|\Psi(c,-\al;xe^{-i\pi})\right|^{-2}.
\eeq
Again the limit $c=0$ leads back to the Laguerre polynomial case.

The explicit form (\ref{EFLP}) simplify in
 \beq
\hspace{-0.3cm}L^\al_{n}(x;c,\mu)\hspace{-0.1cm}
=\hspace{-0.1cm}\frac{(c+\al+1)_n}{n!}
\hspace{-0.1cm}\sum_{k=0}^n\hspace{-0.1cm}
\frac{(-n)_k}{(c+1)_k(c+\al+1)_k}\fx{3}{2}{k-n,c,c+\al+1}
{c+k+1,c+\al+k+1}{1}x^k\hspace{-0.05cm}.
\eeq

The coefficients of the differential equation satisfied by the
zero related Laguerre polynomials are
\beqa\label{ZR1}
c_4 &&\hspace{-0.6cm}=   x^2(2(2 n+1)x+D),\hspace{3.6cm}
c_3 =  2 x(3(2 n+1)x+2 D),\\
c_2 &&\hspace{-0.6cm}=-2(2 n+1)x^3+(8 n (F+1)+8 c+D)x^2\\ \nonumber
&&\hspace{-0.6cm}\hspace{1.75cm}-(4(\al^2-1)n-2  \al^2-
D(4c+1)-1)x-1/4D(D^2-9),\\
c_1 &&\hspace{-0.6cm}=   8(2 n+1)x^2-4(2 n (F+1)+2 c-D)x-2
D(n(D+3)+6 c+2 D),\\
c_0 &&\hspace{-0.6cm}=  n(n+1)(2(2 n+1)x+3 D),
\eeqa
with
\beq\label{ZR6}
F=n+2c+\al,\hspace{1cm}D =  1+2 \al,
\eeq
and are no longer invariant under the transformation ${\cal T}$.

\subsubsection{Limit $\mu=c+\al$ }\label{NLC}
This is a new simple case of CAL polynomials lacking in \cite{ISM88}.
\beqa
{\goth L}_n^\al(x,c)=e^{-x}&&\hspace{-0.6cm}\left\{
\fx{1}{1}{-c}{\al}{x}\fx{1}{1}{-n-c-\al}{1-\al}{x}\right.\nonumber\\
&&\hspace{-0.6cm}+\frac{(c+\al)(c+\al+1)_n}{\al(1-\al)(c+1)_n}x
\left.\fx{1}{1}{1-c-\al}{2-\al}{x}\fx{1}{1}{-n-c}{1+\al}{x}\right\},
\eeqa
and the measure is obtained using \cite[(10) page 258]{HTF1}
\beq
\phi'(x)=\frac{x^{\al-1} e^{-x}}{\Ga(1+c,1+c+\al)}
\left|\Psi(c+1,2-\al;xe^{-i\pi})\right|^{-2}.
\eeq
In this case the limit $c=0$ does not lead back to the Laguerre polynomial case
but to the co-recursive Laguerre one with $\mu=\al$.

The explicit form is
 \beq
\hspace{-0.3cm}L^\al_{n}(x;c,\mu)\hspace{-0.1cm}
=\hspace{-0.1cm}\frac{(c+\al+1)_n}{n!}
\hspace{-0.1cm}\sum_{k=0}^n\hspace{-0.1cm}
\frac{(-n)_k}{(c+1)_k(c+\al+1)_k}\fx{3}{2}{k-n,c+1,c+\al}
{c+k+1,c+\al+k+1}{1}x^k\hspace{-0.05cm}.
\eeq

The coefficients of the differential equation (\ref{CALDE})  satisfied by the
 polynomials ${\goth L}_n^\al(x,c)$ are obtained from
(\ref{ZR1}--\ref{ZR6}) changing only $D=1+2\al$ by $D=1-2\al$.

\section{The case of Jacobi polynomials}     \label{SJC}

We now present some results on the CAJ polynomials.
The recurrence relation of the associated
Jacobi polynomials $P_n^{\al,\be}(x;c) $ is \cite{WIM87}
\beqa
(2n+2c+\al+\be+1)&&\hspace{-0.6cm}
[(2n+2c+\al+\be+2)(2n+2c+\al+\be)x+\al^2-\be^2]p_n
\nonumber\\
&&\hspace{-0.6cm}\hspace{0.6cm}=2(n+c+1)(n+c+\al+\be+1)(2n+2c+\al+\be)p_{n+1}\\
&&\hspace{-0.6cm}\hspace{1.cm}+2(n+c+\al)
(n+c+\be)(2n+2c+\al+\be+2)p_{n-1}.\nonumber
\label{RRAJ}\eeqa
We can note the invariance of the recurrence relation (\ref{RRAJ}) under
the transformation ${\cal T}'$ defined by
\beq
{\cal T}'(c,\al,\be)=(c+\al+\be,-\al,-\be).
\eeq
All polynomials
satisfying  (\ref{RRAJ}), for which the initial conditions are symmetric
in $\al$ and $\be$ and invariant under ${\cal T}'$ have the property
\beq\label{PP1}
p_n^{-\al,-\be}(x,c+\al+\be)=p_n^{\al,\be}(x,c).
\eeq

As in \cite {WIM87} we use the more convenient shifted polynomials defined as
\beq
R_n^{\al,\be}(x;c)=P_n^{\al,\be}(2x-1;c).
\eeq
Due also to the properties of the recurrence relation (\ref{RRAJ}) we have
\beq\label{PP2}
R_n^{\al,\be}(x;c)=(-1)^n R_n^{\be,\al}(1-x;c).
\eeq

\subsection{Explicit representation for the CAJ polynomials}   \label{SSJEF}
A solution of the recurrence relation satisfied by the $R_n^{\al,\be}(x;c)$ in
terms of the
hypergeometric function is \cite[page 280]{LUK1}
\beq
u_n=\frac{(c+\al+1)_n}{(c+1)_n}\fx{2}{1}{-n-c,n+c+\al+\be+1}{1+\al}{1-x},
\eeq
and an other linearly independent solution is given by
\beq
v_n={\cal T}'u_n=\frac{(c+\be+1)_n}{(c+\al+\be+1)_n}
\fx{2}{1}{-n-c-\al-\be,n+c+1}{1-\al}{1-x}.
\eeq
The functions $u_n$ and $x^{-\be}(1-x)^{-\al}v_n$ are two independent solutions
of the second-order differential equation
\beq\label{DEF}
x(1-x)y''(x)+[1+\be-(\al+\be+2)x]y'(x)+(n+c)(n+c+\al+\be+1)y(x)=0.
\eeq

The associated Jacobi polynomials are defined by (\ref{RRAJ}) and the initial
condition
\beq\label{ICJ}
P_{-1}^{\al,\be}(x;c)=0,\hspace{1cm}P_0^{\al,\be}(x;c)=1.
\eeq
This gives for $P_1^{\al,\be}(x;c)$
\beqa
P_1^{\al,\be}(x;c)&=&\frac{(2c+\al+\be+1)(2c+\al+\be+2)}{2(c+1)(c+\al+\be+1)}
\nonumber\\
&&\hspace{5.cm}\times\left[x+\frac{\al^2-\be^2}
{(2c+\al+\be)(2c+\al+\be+2)}\right].
\eeqa
The CAJ polynomials $P_n^{\al,\be}(x;c,\mu)$ satisfy the recurrence relation
(\ref{RRAJ}) with a shift $\mu$ on the first monic polynomial.
This corresponds to the initial condition on the shifted CAJ polynomials
$R_n^{\al,\be}(x;c,\mu)$
\beq\label{ICCJ}
R_{-1}^{\al,\be}(x;c,\mu)=D
=-\frac{(2c+\al+\be)(2c+\al+\be+1)}{2(c+\al)(c+\be)}\mu,\
R_{0}^{\al,\be}(x;c,\mu)=1.
\eeq
If $c+\al\ra 0$ or $c+\be\ra 0$ this initial condition in (\ref{RRAJ}) leads
nevertheless to a shift $\mu$ on the value of $x$ in $P_1^{\al,\be}(x;c)$.

As in section \ref{SSLEF} writing
\beq\label{RN}
R_{n}^{\al,\be}(x;c,\mu)=Au_n+Bv_n
\eeq
and using (\ref{ICCJ}) we obtain
\beq\label{A}
A=\frac{1}{\Delta}\left[D v_0-v_{-1}\right]\hspace{1cm}{\rm and}\hspace{1cm}
B=-\frac{1}{\Delta}\left[D u_0-u_{-1}\right],
\eeq
where $\Delta$ is easily  calculated using the fact that  $u_n$ and
$x^{-\be}(1-x)^{-\al}v_n$ are two independent solutions of (\ref{DEF}),
\beq
\Delta=u_{-1}v_0-u_0v_{-1}=-\frac{\al(2c+\al+\be)}{(c+\al)(c+\be)}.
\eeq
The condition $\De\neq 0$ leads to $\al\neq 0$ and $2c+\al +\be\neq 0$.
We can note the invariance of $D$ under ${\cal T}'$ and that $B={\cal T}'A$.

Grouping the two ${_2F_1}$ involved in the expression (\ref{A}) of $A$ gives
\beq
A=\frac{c+\al+\be-D(c+\be)}{\al(2c+\al+\be)}
\fx{3}{2}{-c-\al-\be,c,F+1}{1-\al,F}{1-x},
\eeq
with
\beq
F=\frac{c[D(c+\be)-c-\al-\be]}{D(c+\be)+c},
\eeq
and the CAJ polynomials could be writed
\beqa
R_{n}^{\al,\be}(x;c,\mu)&&\hspace{-0.6cm}=(1+{\cal T}')
\frac{c+\al+\be-D(c+\be)}{\al(2c+\al+\be)}
(c+\al)\frac{(c+\al+1)_n}{(c+1)_n}\nonumber\\
&&\hspace{-0.6cm}\hspace{-0.9cm}\times\fx{3}{2}{-c-\al-\be,c,F+1}{1-\al,F}{1-x}
\fx{2}{1}{-n-c,n+c+\al+\be+1}{1+\al}{1-x}.
\label{CAJP}\eeqa
We will use this expression of the CAJ polynomials in \ref{FODECAJ} to obtain a
fourth-order differential equation of them.

Transforming the ${_2F_1(1-x)}$ in (\ref{RN}) by \cite[Eq. 1, page 108]{HTF1}
one obtains
with little algebra
\beqa
\hspace{1.cm}R_{n}^{\al,\be}(x;c,\mu)&&\hspace{-0.6cm}=(1+{\cal T}')
\frac{(-1)^nc(c+\al)}{\be(2c+\al+\be)}\nonumber\\
&&\hspace{-0.6cm}\times\frac{(c+\al+1)_n}{(c+\al+\be+1)_n}
\fx{2}{1}{-n-c-\al-\be,n+c+1}{1-\be}{x}\label{CAJW}\\
&&\hspace{-0.6cm}\times\left[\frac{c+\be}{c}D
\fx{2}{1}{-c,c+\al+\be+1}{1+\be}{x}
-\fx{2}{1}{1-c,c+\al+\be}{1+\be}{x}\right].\nonumber
\eeqa
This formula generalizes the one of \cite[Eq. 28]{WIM87} to the case of the CAJ
polynomials.
As the explicit form of the CAL polynomials (\ref{CALP}) the representation
(\ref{CAJP})
and  (\ref{CAJW}) are valid
only for $\al\neq 0,\pm 1,\pm 2\ldots$ and $\be\neq 0,\pm 1,\pm 2\ldots$
but can be extended by limiting processes.

We obtain an explicit formula following the same way as in \cite{WIM87}. We
first use
\cite[Eq. 14, page 87]{HTF1} for each product of ${_2F_1}$ in (\ref{CAJW}) to
obtain four
series involving gamma functions and a ${_4F_3}$. For two of them we use
\cite[Eq. 1, page 56]{BAI72}. The next step is to use for each ${_4F_3}$ twice
\cite[Eq. 3, page 62]{BAI72}. After numerous cancellations only two series of
${_4F_3}$
remains we can group to obtain the following explicit form
 \beqa
R_{n}^{\al,\be}(x;c,\mu)=(-1)^n\frac{(2c+\al+\be+1)_n(\be+c+1)_n}
{n!(c+\al+\be+1)_n}
\sum_{k=0}^n&&\hspace{-0.6cm}
\frac{(-n)_k(n+2c+\al+\be+1)_k}{(c+1)_k(c+\be+1)_k}\nonumber\\
&&\hspace{-0.6cm}\hspace{-5.4cm}\times
\fx{5}{4}{k-n,n+k+2c+\al+\be+1,c,c+\be,G+1}
{c+k+1,c+\be+k+1,2c+\al+\be+1,G}{1}x^k,\label{EFCAJ}
\eeqa
where
\beq\label{G}
G=\frac{2c(c+\be)(2c+\al+\be)}{2c(c+\be)+\mu(2c+\al+\be)(2c+\al+\be+1)}.
\eeq

\subsection{Generating function}\label{SSCAJGF}
One can obtain a generating function of the CAJ polynomials following the same
strategy as
in \cite{WIM87} for the associated one. Let ${\cal G}(x,w)$ be a generating
function of
$R_{n}^{\al,\be}(x;c,\mu)$
\beq
{\cal G}(x,w)=\sum_{n=0}^\nf\frac{(c+1)_n(c+\al+\be+1)_n}{n!(2c+\al+\be+2)_n}
w^nR_{n}^{\al,\be}(x;c,\mu).
\eeq
Starting from the form (\ref{CAJP}) for the $R_{n}^{\al,\be}(x;c,\mu)$ it
follows
\beqa
\hspace{-0.6 cm}{\cal G}(x,w)&&\hspace{-0.6cm}=(1+{\cal T}')
\frac{c+\al+\be-D(c+\be)}{\al(2c+\al+\be)}
(c+\al)\fx{3}{2}{-c-\al-\be,c,F+1}{1-\al,F}{1-x}\nonumber\\
&&\hspace{-0.6cm}\times\sum_{n=0}^\nf
\frac{(c+\al+1)_n(c+\al+\be+1)_n}{n!(2c+\al+\be+2)_n}
w^n\fx{2}{1}{-n-c,n+c+\al+\be+1}{1+\al}{1-x},
\eeqa
and using \cite[Th. 4]{WIM87} we obtain
\beqa
{\cal G}(x,w)&&\hspace{-0.6cm}=(1+{\cal
T}')\frac{c+\al+\be-D(c+\be)}{\al(2c+\al+\be)}
\nonumber\\&&\hspace{-0.6cm}\times(c+\al)
\left[\frac{2}{w(Z_2+1)}\right]^{c+\al+\be+1}
\fx{3}{2}{-c-\al-\be,c,F+1}{1-\al,F}{1-x}
\\&&\hspace{-0.6cm}\times
\fx{2}{1}{-c,c+\al+\be+1}{1+\al}{\frac{1-Z_1}{2}}
\fx{2}{1}{c+\al+1,c+\al+\be+1}{2c+\al+\be+2}{\frac{2}{1+Z_2}},\nonumber
\eeqa
where
\beq
Z_1=\frac{1-\sqrt{(1+w)^2-4wx}}{w},\hspace{1cm}Z_2
=\frac{1+\sqrt{(1+w)^2-4wx}}{w},
\eeq
which generalize the already exotic generating function \cite[Eq. 75]{WIM87}.

\subsection{Spectral measure}\label{SSCAJM}
The Stieltjes transform of the measure of the shifted associated Jacobi
polynomials
$R_{n}^{\al,\be}(x;c)$ is \cite[Eq. 63--64]{WIM87}
\beq\label{SMAJ}
s(p)=\frac{1}{p}\,\fx{2}{1}{c+1,c+\be+1}{2c+\al+\be+2}{\frac{1}{p}}
{\fx{2}{1}{c,c+\be}{2c+\al+\be}{\frac{1}{p}}}^{-1}.
\eeq

The CAJ polynomials $R_{n}^{\al,\be}(x;c,\mu)$ satisfy the same recurrence
relations
as the $R_{n}^{\al,\be}(x;c)$ with a shift $\mu$ on the first monic polynomials
\beq
R_{1}^{\al,\be}(x;c,\mu)-R_{1}^{\al,\be}(x;c)
=\frac{(2c+\al+\be+1)(2c+\al+\be+2)}
{2(c+1)(c+\al+\be+1)}\mu.
\eeq

Using continued J-fractions \cite{ISM90,HTF2,SHO50} whose denominators are
$R_{n}^{\al,\be}(x;c,\mu)$
and $R_{n}^{\al,\be}(x;c)$ we can derive for the Stieltjes transform of the
measure
of the CAJ polynomials
\beqa
s(p;\mu)&&\hspace{-0.6cm}=s(p)\left(1+\frac{\mu}{2}s(p)\right)^{-1}\nonumber\\
&&\hspace{-0.6cm}=\frac{\fx{2}{1}{c+1,c+\be+1}{2c+\al+\be+2}
{\displaystyle\frac{1}{p}}}
{p\,\fx{2}{1}{c,c+\be}{2c+\al+\be}{\displaystyle\frac{1}{p}}
+{\displaystyle\frac{\mu}{2}}\,\fx{2}{1}{c+1,c+\be+1}{2c+\al+\be+2}
{\displaystyle\frac{1}{p}}},
\label{STCAJ}
\eeqa
that we can also write using contiguous relations
\beqa
s(p;\mu)&&\hspace{-0.6cm}=\fx{2}{1}{c+1,c+\be+1}{2c+\al+\be+2}
{\displaystyle\frac{1}{p}}
\nonumber\\
&&\hspace{-0.6cm}\times\label{MUP}
\left[\left(\frac{c+\al}{2c+\al+\be}-\frac{2c+\al+\be+1}{2c}\mu\right)
\fx{2}{1}{c,c+\be}{2c+\al+\be+1}{\displaystyle\frac{1}{p}}\right.\\
&&\hspace{-0.6cm}\hspace{0.6cm}\left.+\left(\frac{c+\be}{2c+\al+\be}
+\frac{2c+\al+\be+1}{2c}\mu\right)
\fx{2}{1}{c,c+\be+1}{2c+\al+\be+1}{\displaystyle\frac{1}{p}}
\right]^{-1}.\nonumber
\eeqa
Grouping the ${_2F_1}$ in (\ref{MUP}) gives the compact formula
\beq\label{STCOM}
s(p;\mu)=\fx{2}{1}{c+1,c+\be+1}{2c+\al+\be+2}{\displaystyle\frac{1}{p}}
\left[\fx{3}{2}{c,c+\be,G+1}{2c+\al+\be+1,G}{\frac{1}{p}}\right]^{-1}.
\eeq
where $G$ is given by (\ref{G}). A sufficient condition for the positivity of
the
denominator in (\ref{STCOM}) on $(1, \nf)$ is
\beq\label{COND}
c\geq 0,\hspace{1cm}c>-\be,\hspace{1cm}\al>-1,\hspace{1cm}
\mu\geq - \frac{ 2 c(c+\be)}{(2 c+\al+ \be)(2 c+\al+ \be+1)},
\eeq
but other conditions are possible.

To obtain the absolutely continuous part of the spectral measure we need to
evaluate
$s^+(p;\mu)-s^-(p;\mu)$ where $s^\pm$
are the values of $s$ above and below the cut [0,1]. Using the analytic
continuation
\cite[Eq. 2, page 108]{HTF1} for each ${_2F_1}$ in (\ref{STCAJ}) we find for
the spectral
measure of the CAJ polynomials
\beqa
\phi'(x)&&\hspace{-0.6cm}=(1-x)^\al x^{\be+2c}
\left| \fx{2}{1}{c,c+\be}{2c+\al+\be}{\frac{e^{i\pi}}{x}}
+\frac{\mu}{2x}\,\fx{2}{1}{c+1,c+\be+1}{2c+\al+\be+2}
{\frac{e^{i\pi}}{x}}\right|^{-2}
\nonumber\\
&&\hspace{-0.6cm}=(1-x)^\al x^{\be+2c}
\left| \fx{3}{2}{c,c+\be,G+1}{2c+\al+\be+1,G}{\frac{e^{i\pi}}{x}}\right|^{-2},
\eeqa
valid at least under the conditions (\ref{COND}).

\subsection{Fourth-order differential equation}\label{FODECAJ}
The method used to  obtain the differential equation satisfied by the
$R_{n}^{\al,\be}(x;c,\mu)$ is the same as in \ref{SSLDE}.
In (\ref{CAJP}) the hypergeometric function ${_2F_1}$  is solution of the
equation
(\ref{DEF})  and the ${_3F_2}$   is of the form
$\fx{3}{2}{a,b,e+1}{d,e}{x}$ which is also solution of the second-order
differential
equation
\beqa
&&\hspace{-0.6cm}x(x-1)\left[(a-e)(b-e)x+e(d-e-1)\right]y''(x)\nonumber\\
&&\hspace{-0.6cm}+\left\{(a-e)(b-e)(a+b+1)x^2
+\left[e(a+b+1)(2d-e-2)-d(ab+e^2)+ab\right]x
\right.\\
&&\hspace{-0.cm}\left.\vphantom{\left[x^2+\right]}
+de(e-d+1)\right\}y'(x)+ab\left[(a-e)(b-e)x+(e+1)(d-e-1)\right]y(x)=0.\nonumber
\eeqa

We don't write here the fourth-order differential equation hardly obtained by
symbolic MAPLE computation.
The coefficients are at most of degree eight in $x$ and it would take
several pages to write them. We give the results only in the following limiting
cases.

\subsection{Particular cases}\label{SSJPC}
\subsubsection{Laguerre case limit}
The limit giving the CAL polynomial case is obtained by the replacement
\beq\label{JTOL}
\left.\begin{array}{ccl}
x&\ra&1-\frac{2x}{\be}\\
\mu&\ra&+\frac{2\mu}{\be}\end{array}\right\}\hspace{1cm}\be\ra\nf,
\eeq
in $P_{n}^{\al,\be}(x;c,\mu)$.
 The representation (\ref{CAJP}) is the more suitable to obtain the form
of the CAL polynomials (\ref{CALP}) using
the Kummer's transformation (\ref{K1F1}) for one of the confluent
hypergeometric
functions and his generalization
\beq\label{K2F2}
\fx{2}{2}{a,e+1}{c,e}{x}=e^x\fx{2}{2}{c-a-1,
\frac{e(c-a-1)}{e-a}+1}{c,\frac{e(c-a-1)}{e-a}}{-x}
\eeq
 for one of the ${_2F_2}$. Note that (\ref{K2F2}) gives (\ref{K1F1}) in the
limit $e\ra\nf$.

\subsubsection{Limit $c=0$}   \label{SSCJP}
In this limit we obtain the co-recursive Jacobi polynomials.
 An explicit form is
\beqa\label{EFCJP}
\hspace{3cm}R_{n}^{\al,\be}(x;\mu)&&\hspace{-0.6cm}=
(-1)^n\frac{(\be+1)_n}{n!}\sum_{k=0}^n
\frac{(-n)_k(n+\al+\be+1)_k}{k!(\be+1)_k}x^k\\
&&\hspace{-0.6cm}\hspace{-4.9cm}\times
\left\{1+\frac{\mu(k-n)(n+k+\al+\be+1)}{2(k+1)(\be+k+1)}
\fx{4}{3}{k-n+1,n+k+\al+\be+2,\be+1,1}
{k+2,\be+k+2,\al+\be+2}{1}\right\},\nonumber
\eeqa
and the spectral measure is given by
\beq\label{SMCJP}
\phi'(x)=(1-x)^\al x^{\be}
\left| 1
+\frac{\mu}{2x}\fx{2}{1}{1,\be+1}{\al+\be+2}{\frac{e^{i\pi}}{x}}\right|^{-2},
\eeq
The limit $\mu=0$ leads back to the Jacobi polynomials.

The fourth-order  differential equation  satisfied by the
co-recursive Laguerre polynomials can be factorized in the limit $c=0$ to
obtain as in \cite{RON89} the factorized (2+2) differential equation
\beqa
\hspace{-2cm}&&\hspace{-0.6cm}\hspace{1.5cm}0=\left[(1-x^2)A(x){\rm
D}^2+\left\{(\be-\al-(\al+\be+4)x)A(x)
-(1-x^2)B(x)\right\}{\rm D}\right.\nonumber\\
&&\hspace{-1cm}\left.+\left\{\vphantom{{\rm
D}^2}n(n+\al+\be+1)-(\al+\be+2)\right\}A(x)
+\left\{\vphantom{{\rm
D}^2}\be-\al-(\al+\be+2)x\right\}B(x)+C(x)\right]\label{FECJ}\\
&&\hspace{-1cm}\hspace{-0.cm}\times\left[(1-x^2){\rm D}^2+\{\vphantom{{\rm
D}^2}
(\al+\be-2)x+\al-\be\}{\rm
D}+n(n+\al+\be+1)+\al+\be\right]R_{n}^{\al,\be}(x;\mu),
\nonumber
\eeqa
where
\beqa
A(x)&&\hspace{-0.6cm}=2(\al+\be)^2(2n+1)(n+\al
+\be+1/2)x^2+2(\al+\be)(4n(n+\al+\be+1)
\nonumber\\&&\hspace{-0.6cm}\times
(-\mu(1+\al+\be)+\al-\be)+(1+\al+\be)(-\mu(\al+\be+2)+2\al-2\be))x\nonumber\\
&&\hspace{-0.6cm}+4n(n+\al+\be+1)(-\mu(1+\al+\be)+\al-\be)^2-(\al+\be)(-2\mu(
1+\al+\be)\nonumber\\
&&\hspace{-0.6cm}\times (\be-\al)-2(\be-\al)^2-3\al-3\be)\nonumber\\
B(x)&&\hspace{-0.6cm}=-(\al+\be)((\al+\be)(8n(n+\al+\be
+1)+3\al+3\be)x+8n(n+\al+\be+1)\nonumber\\&&\hspace{-0.6cm}\times
(-\mu(1+\al+\be)+\al-\be)-2(\al+\be+2)(1+\al+\be)\mu
+(\al-\be)(3\al+3\be+4))\nonumber\\
C(x)&&\hspace{-0.6cm}=-(\al+\be)((\al+\be)(\al+\be+2)
(2n(n+\al+\be+1)+\al+\be-1)x^2\nonumber\\
&&\hspace{-0.6cm}+2(n(n+\al+\be+1)(-\mu(\al+\be+1)(\al+\be-4)
+2(\al-\be)(\al+\be-2))\nonumber\\
&&\hspace{-0.6cm}+(\al+\be)(\al-\be)(\al+\be-1))x
-2n(n+\al+\be+1)(-\mu(\al+\be+1)(\be-\al)
\nonumber\\&&\hspace{-0.6cm}-(\al-\be)^2+6\al+6\be)
+(\al+\be)((\al-\be)^2-3\al-3\be-6))\nonumber
\eeqa

\subsubsection{Limit $c=-\al-\be$}
Due to the $\cal T'$ invariance of (\ref{RRAJ}) we obtain in this limit the
special
case of CAJ polynomials for which
\beq
R_{n}^{\al,\be}(x;-\al-\be,\mu)=R_{n}^{-\al,-\be}(x;\mu).
\eeq
All the results are obtained from \ref{SSCJP} by changing $\al$ to $-\al$ and
$\be$
to $-\be$.

\subsubsection{Limit $c=-\be$}\label{CBE}
The explicit form (\ref{EFCAJ}) simplifies in the same way as in the case
$c=0$. One obtain
\beqa
\hspace{2.cm}R_{n}^{\al,\be}(x;-\be,\mu)&&\hspace{-0.6cm}=
(-1)^n\frac{(\al-\be+1)_n}{(\al+1)_n}\sum_{k=0}^n
\frac{(-n)_k(n+\al-\be+1)_k}{k!(1-\be)_k}x^k\\
&&\hspace{-0.6cm}\hspace{-4.6cm}\times
\left\{1+\frac{\mu(k-n)(n+k+\al-\be+1)}{2(k+1)(1-\be+k)}
\fx{4}{3}{k-n+1,n+k+\al-\be+2,1-\be,1}
{k+2,2-\be+k,\al-\be+2}{1}\right\}.\nonumber
\eeqa
Comparing this form with the explicit form of the co-recursive Jacobi
polynomials
(\ref{EFCJP}) one see
\beq
R_{n}^{\al,\be}(x;-\be,\mu)=\frac{n!(\al-\be+1)_n}{(\al+1)_n(1-\be)_n}
R_{n}^{\al,-\be}(x;\mu),
\eeq
of course the spectral measure and  the fourth-order  differential equation are
obtained
from (\ref{EFCJP}) and (\ref{FECJ}) changing $\be$ to $-\be$.

\subsubsection{Limit $c=-\al$}
This case is the $\cal T'$ transform of the preceding case.
All the results are obtained from \ref{CBE} by changing $\al$ to $-\al$ and
$\be$
to $-\be$.

\subsubsection{Limit $\mu=0$}

In this limit we obtain the associated Jacobi polynomials studied in
\cite{WIM87}.
The form \cite[Eq. 28]{WIM87} is obtain directly using (\ref{CAJW}) but
 a  form slightly different is
\beqa
R_{n}^{\al,\be}(x;c)&&\hspace{-0.6cm}=(1+{\cal T}')
\frac{(c+\al)(c+\al+\be)}{\al(2c+\al+\be)}
\frac{(c+\al+1)_n}{(c+1)_n}\nonumber\\
&&\hspace{-0.6cm}\hspace{-0.8cm}\times\fx{2}{1}{1-c-\al-\be,c}{1-\al}{1-x}
\fx{2}{1}{-n-c,n+c+\al+\be+1}{1+\al}{1-x},
\label{AJP}\eeqa
The explicit form \cite[Eq. 19]{WIM87} is easily obtain starting from
(\ref{EFCAJ}) with
$G=2c+\al+\be$, the ${_5F_4}$ reducing to a ${_4F_3}$.
Obviously the limit $c=0$ leads back to the Jacobi polynomials.

The coefficients of the differential equation (\ref{CALDE}) satisfied by the
associated Jacobi polynomials are
\beqa
c_4 &&\hspace{-0.6cm}=   x^2(x-1)^2,\hspace{1.6cm}c_3 =  5 x(x-1)(2 x-1),\\
c_2&&\hspace{-0.6cm}=\left(24-(n+1)^2-A\right)x(x-1)-B x-\be^2+4,\\
c_1 &&\hspace{-0.6cm}=-3/2\left(\vphantom{(C^2)}(3A+(n+3)(n-1))(2
x-1)+B\right),\\
c_0 &&\hspace{-0.6cm}=  n(n+2)A,
\eeqa
with
\beq\label{C}
A =(C+n+1)(C+n-1),\hspace{0.5cm}
B =(\al-\be)(\al+\be),\hspace{0.5cm}
C = 2 c+\al+ \be.
\eeq
This result was also first given by Hahn \cite[Eq. 20]{HAH40}.
Note the $\cal T'$ invariance of $A,\ B$ and $C$ leading to the invariance of
the $c_i$
 more obvious than in \cite[Eq. 47--48]{WIM87}.

\subsubsection{Limit $\mu=\frac{2c(c+\al)}{(2c+\be+\al)(2c+\be+\al+1)}$}
\label{ZRJ}
In this limit the symmetry ${\cal T}'$ is broken.
We obtain the zero-related Jacobi polynomials studied in \cite{ISM91}.
An explicit form is
 \beqa
{\cal
R}_{n}^{\al,\be}(x;c)=(-1)^n
\frac{(2c+\al+\be+1)_n(\be+c+1)_n}{n!(c+\al+\be+1)_n}
\sum_{k=0}^n&&\hspace{-0.6cm}
\frac{(-n)_k(n+2c+\al+\be+1)_k}{(c+1)_k(c+\be+1)_k}\nonumber\\
&&\hspace{-0.6cm}\hspace{-6.4cm}\times\fx{4}{3}{k-n,n+k+2c+\al+\be+1,c,c+\be+1}
{c+k+1,c+\be+k+1,2c+\al+\be+1}{1}x^k.
\eeqa
The limit $c=0$ leads back to the Jacobi polynomials and the limit defined in
(\ref{JTOL})
gives the zero-related Laguerre polynomials (\ref{ZRL}), using Kummer's
transformations.
The spectral measure is
\beq
\phi'(x)=(1-x)^\al x^{\be+2c}
\left| \fx{2}{1}{c,c+\be+1}{2c+\al+\be+1}{\frac{e^{i\pi}}{x}}\right|^{-2}.
\eeq

The coefficients of the differential equation (\ref{CALDE})
satisfied by the polynomials ${\cal R}_{n}^{\al,\be}(x;c)$ are
\beqa
c_4 =&&\hspace{-0.6cm}   x^2(x-1)^2(Ax+D),\label{Z1}\\
c_3= &&\hspace{-0.6cm} x(x-1)\left(8Ax^2-3(A-3D)x-4D\right),\\
c_2 =&&\hspace{-0.6cm}-1/2A(A+2C^2-29)x^3
+\left(1/2A(A+2C^2-2B-23)-D(C^2-19)\right)x^2\\
&&\hspace{-0.6cm}-1/4\left(A(D+1)(D-3)-2D(2C^2-2B+D-35)\right)x
-1/4D(D^2-9),\nonumber\\
c_1 =&&\hspace{-0.6cm}-A(A+2C^2-5)x^2
+1/4\left(A(A+2C^2-2B-5D-5)-3D(4C^2-11)\right)x\\
&&\hspace{-0.6cm}+1/4D\left((D+3)A+6C^2-6B+3D-15\right),\nonumber\\
c_0 = &&\hspace{-0.6cm} 2n(n+1)(C+n)(C+n+1)(Ax+3D),
\eeqa
where $B$ and $C$ are defined in (\ref{C}) and
\beq\label{Z6}
A =(2 n+1)(1+2 C+2 n),\hspace{1cm}D =  1+2 \be.
\eeq

\subsubsection{Limit
$\mu=\frac{2(c+\be)(c+\al+\be)}{(\be+\al+2c)(\be+\al+2c+1)}$}
This case is the ${\cal T}'$ transform of the case \ref{ZRJ}. The explicit form
is
 \beqa
\hspace{0.8cm}{\goth R}_{n}^{\al,\be}(x;c)
=(-1)^n\frac{(2c+\al+\be+1)_n(\al+c+1)_n}{n!(c+1)_n}
\sum_{k=0}^n&&\hspace{-0.6cm}
\frac{(-n)_k(n+2c+\al+\be+1)_k}{(c+\al+\be+1)_k(c+\al+1)_k}\nonumber\\
&&\hspace{-0.6cm}\hspace{-7.4cm}\times
\fx{4}{3}{k-n,n+k+2c+\al+\be+1,c+\al+\be,c+\al+1}
{c+\al+\be+k+1,c+\al+k+1,2c+\al+\be+1}{1}x^k.
\eeqa
The limit (\ref{JTOL}) leads back to the Laguerre case \ref{NLC} and the limit
$c=0$ to the
co-recursive Jacobi polynomials with $\mu=\frac{2\be}{\al+\be+1}$.
The spectral measure is
\beq
\phi'(x)=(1-x)^{\al} x^{\be+2c}
\left| \fx{2}{1}{c+1,c+\be}{2c+\al+\be+1}{\frac{e^{i\pi}}{x}}\right|^{-2}.
\eeq

The coefficients
of the differential equation (\ref{CALDE}) satisfied by the polynomials
${\goth R}_{n}^{\al,\be}(x;c)$ are obtained from (\ref{Z1}--\ref{Z6}) changing
only
$D =  1+2 \be$ by $D =  1-2 \be$.

\subsubsection{Limit $\mu = - \frac{ 2 c(c+\be)}{(2 c+\al+ \be)(2 c+\al+
\be+1)}$}
\label{NRJ}
The symmetry ${\cal T}'$ is also broken. We obtain a new simple
case of CAJ polynomials.
An explicit form is
 \beqa
\widetilde{{\cal R}}_{n}^{\al,\be}(x;c)
=(-1)^n\frac{(2c+\al+\be+1)_n(\be+c+1)_n}{n!(c+\al+\be+1)_n}
\sum_{k=0}^n&&\hspace{-0.6cm}
\frac{(-n)_k(n+2c+\al+\be+1)_k}{(c+1)_k(c+\be+1)_k}\nonumber\\
&&\hspace{-0.6cm}\hspace{-6.4cm}\times\fx{4}{3}{k-n,n+k+2c+\al+\be+1,c,c+\be}
{c+k+1,c+\be+k+1,2c+\al+\be+1}{1}x^k.
\eeqa
and the spectral measure
\beq
\phi'(x)=(1-x)^\al x^{\be+2c}
\left| \fx{2}{1}{c,c+\be}{2c+\al+\be+1}{\frac{e^{i\pi}}{x}}\right|^{-2}.
\eeq

The coefficients of the differential equation (\ref{CALDE})
satisfied by the polynomials $\widetilde{{\cal R}}_{n}^{\al,\be}(x;c)$ are
\beqa
c_4 =&&\hspace{-0.6cm}   x^2(x-1)^2\left(A(x-1)-D\right),\label{ZZ1}\\
c_3= &&\hspace{-0.6cm}  x(x-1)\left(8 A  x^2-(13 A+9 D)x+5 A+5 D\right),\\
c_2 =&&\hspace{-0.6cm}-1/2 A(A+2  C^2-29)x^3
+\left(A(A+2  C^2-B-32)+D(C^2-19)\right)x^2\\
&&\hspace{-0.6cm}-1/2\left(A(A+2  C^2+2 \be^2-2 B-43)+D(2  C^2-2 B-D-41
)\right)x \nonumber\\
&&\hspace{-0.6cm} +(\be^2-4)(A+D),\nonumber\\
c_1 =&&\hspace{-0.6cm}-2 A(A+2  C^2-5)x^2 \\
&&\hspace{-0.6cm}+  1/2\left( A(7 A+14  C^2-2 B+ 5 D-35) +D (3
C^2-33)\right)x\nonumber\\
&&\hspace{-0.6cm}-A(3/2 A+3  C^2-B-2  \al^2-7)-3 D(C^2-B-1/2D-3),\nonumber\\
c_0 = &&\hspace{-0.6cm}  n(n+1)(C+n)(C+n+1)(A(x-1)-3 D),
\eeqa
where $B$ and $C$ are defined in (\ref{C}) and
\beq\label{ZZ6}
A =(2 n+1)(1+2 C+2 n),\hspace{1cm}D =  1+2 \al.
\eeq

\subsubsection{Limit $\mu = -
\frac{2( c+ \al)(c+\al+\be)}{(2 c+\al+ \be)(2 c+\al+ \be+1)}$}
This case is the ${\cal T}'$ transform of the case \ref{NRJ}. The explicit form
is
 \beqa
\widetilde{{\goth R}}_{n}^{\al,\be}(x;c)
=(-1)^n\frac{(2c+\al+\be+1)_n(\al+c+1)_n}{n!(c+1)_n}
\sum_{k=0}^n&&\hspace{-0.6cm}
\frac{(-n)_k(n+2c+\al+\be+1)_k}{(c+\al+\be+1)_k(c+\al+1)_k}\nonumber\\
&&\hspace{-0.6cm}\hspace{-6.4cm}\times
\fx{4}{3}{k-n,n+k+2c+\al+\be+1,c+\al+\be,c+\al}
{c+\al+\be+k+1,c+\al+k+1,2c+\al+\be+1}{1}x^k.
\eeqa
and the spectral measure
\beq
\phi'(x)=(1-x)^{\al-2} x^{\be+2c+2}
\left| \fx{2}{1}{c+1,c+\be+1}{2c+\al+\be+1}{\frac{e^{i\pi}}{x}}\right|^{-2}.
\eeq

The coefficients
of the differential equation (\ref{CALDE}) satisfied by the polynomials
$\widetilde{{\goth R}}_{n}^{\al,\be}(x;c)$ are obtained from
(\ref{ZZ1}--\ref{ZZ6})
by changing only
$D =  1+2 \al$ by $D =  1-2 \al$.

\subsection{Conclusion}\label{SCONC}
We end by  brief remarks. In this article we have studied properties of the
co-recursive
associated Laguerre and Jacobi polynomials which are of interest in the
resolution of some
birth and death processes with and without absorption. For few values of the
co-recursivity
parameter we obtain polynomial families for which the results are of the same
complexity as
the corresponding associated polynomials. For the CAL polynomials we find the
two expected
cases corresponding to $\mu=\mu_0$ (zero-related polynomials) and the new dual
case $\mu=\la_{-1}$ \cite{LET86}. For the CAJ polynomials,
due to the properties (\ref{PP1}) and (\ref{PP2})
we have two more cases corresponding to $\mu=-{\cal T''}\mu_0$ and $\mu=-{\cal
T''}\la_{-1}$
where the transformation ${\cal T''}$ is defined by
\beq
{\cal T''}(c,\al,\be)=(c+\al+\be,-\be,-\al).
\eeq

 In some cases the fourth-order differential equations satisfied by the
polynomials
studied above are factorizable (co-recursive and associated of order one) but
we
don't find factorization either for the co-recursive associated polynomials or
for the
associated one.
 Of course this is not a proof that the conjectures on this factorizability
made in
\cite{RON89} are wrong.

\vspace{1cm}
\noindent{\bf Acknowledgement.} We thank Galliano Valent and Pascal Maroni for
stimulating
 discussions during the preparation of this article.

\bibliographystyle{plain}

\begin{thebibliography}{10}

\bibitem{AND85}
G.~E. Andrews and R.~Askey.
\newblock Classical orthogonal polynomials.
\newblock In C.~Brezinski et~al, editor, {\em Polyn\^{o}mes Orthogonaux et
  Applications}, volume 1171, pages 36--62. Springer Verlag. Berlin, 1985.

\bibitem{ASK84B}
R.~Askey and J.~Wimp.
\newblock Associated {Laguerre} and {Hermite} polynomials.
\newblock {\em Proc. Royal Soc. Edinburgh}, 96A:$\ $15--37, 1984.

\bibitem{ATK81}
F.~V. Atkinson and W.~N. Everitt.
\newblock Orthogonal polynomials which satisfy second-order differential
  equations.
\newblock In P.~Butzer and F.~Feher, editors, {\em Christoffel Festscrift},
  Basel, 1981. {Birkha\"{u}ser-Verlag}.

\bibitem{BAI72}
W.~N. Bailey.
\newblock {\em Generalized Hypergeometric Series}.
\newblock Cambridge University Press, reprinted by Hafner, New York, 1972.

\bibitem{BEL91}
S.~Belmehdi and A.~Ronveaux.
\newblock The fourth-order differential equation satisfied by the associated
  orthogonal polynomials.
\newblock {\em Rendiconti di Matematica Roma (Serie VII)}, 11:$\ $313--326,
  1991.

\bibitem{MAP88}
B.~W. Char, K.~O. Geddes, G.~H. Gonnet, M.~B. Monagan, and S.~M. Watt.
\newblock {\em MAPLE reference manual}.
\newblock WATCOM publications limited, University of Waterloo, Canada, 1988.

\bibitem{HTF1}
A.~{Erd\' elyi}, W.~Magnus, F.~Oberhettinger, and F.~G. T{ricomi}.
\newblock {\em Higher Transcendental Functions}, volume~I.
\newblock McGraw-Hill, New York, 1953.

\bibitem{HTF2}
A.~{Erd\' elyi}, W.~Magnus, F.~Oberhettinger, and F.~G. Tr{icomi}.
\newblock {\em Higher Transcendental Functions}, volume~II.
\newblock McGraw-Hill, New York, 1953.

\bibitem{HAH40}
W.~Hahn.
\newblock {\em { \"Uber} {Orthogonalpolynome mit drei Parametern}}, volume~5.
\newblock Deutche. Math., 1940--41.

\bibitem{ISM79}
M.~E.~H. Ismail and D.~H. Kelker.
\newblock Special functions, {Stieltjes} transforms and infinite divisibility.
\newblock {\em SIAM J. Math. Anal.}, 10:$\ $884--901, 1979.

\bibitem{ISM90C}
M.~E.~H. Ismail, J.~Letessier, D.~Masson, and G.~Valent.
\newblock Birth and death processes and orthogonal polynomials.
\newblock In P.~Nevai, editor, {\em Orthogonal Polynomials: Theory and
  Practice}, volume 294, pages 229--255. NATO ASI series C, 1990.

\bibitem{ISM88}
M.~E.~H. Ismail, J.~Letessier, and G.~Valent.
\newblock Linear birth and death models and associated {Laguerre} polynomials.
\newblock {\em J. Approximation Theory}, 56:$\ $337--348, 1988.

\bibitem{ISM89}
M.~E.~H. Ismail, J.~Letessier, and G.~Valent.
\newblock Quadratic birth and death processes and associated continuous dual
  {Hahn} polynomials.
\newblock {\em SIAM J. Math. Anal.}, 20:$\ $727--737, 1989.

\bibitem{ISM90}
M.~E.~H. Ismail, J.~Letessier, G.~Valent, and J.~Wimp.
\newblock Two families of associated {Wilson} polynomials.
\newblock {\em Canadian J. Math.}, 42:$\ $659--695, 1990.

\bibitem{ISM91}
M.~E.~H. Ismail and D.~R. Masson.
\newblock Two families of orthogonal polynomials related to {Jacobi}
  polynomials.
\newblock {\em Rocky Mountain J. Math.}, 21:$\ $?--?, 1991.

\bibitem{KAR58}
S.~Karlin and J.~McGregor.
\newblock The differential equations of birth and death processes and the
  {Stieltjes} moment problem.
\newblock {\em Transactions Amer. Math. Soc.}, 85:$\ $489--546, 1958.

\bibitem{KAR82}
S.~Karlin and S.~{Tavar\'e}.
\newblock A diffusion process with killing: the time to formation of recurrent
  deleterious mutant genes.
\newblock {\em Stochastic Processes and their Applications}, 13:$\ $249--261,
  1982.

\bibitem{KAR82B}
S.~Karlin and S.~{Tavar\'e}.
\newblock Linear birth and death processes with killing.
\newblock {\em J. Appl. Prob.}, 19:$\ $477--487, 1982.

\bibitem{LAB85}
J.~Labelle.
\newblock Tableau d'{Askey}.
\newblock In C.~Brezinski et~al, editor, {\em Polyn\^{o}mes Orthogonaux et
  Applications}, volume 1171, page XXXVI. Springer Verlag. Berlin, 1985.

\bibitem{LET86}
J.~Letessier and G.~Valent.
\newblock Dual birth and death processes and orthogonal polynomials.
\newblock {\em SIAM J. Appl. Math.}, 46:$\ $393--405, 1986.

\bibitem{LUK1}
Y.~L. {Luke}.
\newblock {\em The special functions and their approximations}, volume~I.
\newblock Academic Press, New York, 1969.

\bibitem{MAR91}
P.~Maroni.
\newblock Une th\'eorie alg\'ebrique des polyn\^{o}mes orthogonaux.
  {Application} aux polyn\^{o}mes orthogonaux semi-classique.
\newblock In C.~Brezinski et~al., editor, {\em Orthogonal polynomials and their
  applications, {\rm IMACS volume 9}}, pages $\ $95--130, 1991.

\bibitem{ORR00}
W.~McF. Orr.
\newblock On the product {$J_m(x)J_n(x)$}.
\newblock {\em Proc. Camb. Phil. Soc.}, 10:$\ $93--100, 1900.

\bibitem{REU57}
G.~E.~H. Reuter.
\newblock Denumerable {Markov} processes and associated semigroups on $l$.
\newblock {\em Acta Math.}, 97:$\ $1--46, 1957.

\bibitem{RON88}
A.~Ronveaux.
\newblock Fourth-order differential equations for numerator polynomials.
\newblock {\em J. Phys. A: Math. Gen.}, 21:$\ $L749--L753, 1988.

\bibitem{RON91}
A.~Ronveaux.
\newblock 4th-order differential equations and orthogonal polynomials of the
  {Laguerre-Hahn} class.
\newblock In C.~Brezinski et~al., editor, {\em Orthogonal polynomials and their
  applications, {\rm IMACS volume 9}}, pages $\ $379--385, 1991.

\bibitem{RON89}
A.~Ronveaux and F.~Marcellan.
\newblock Co-recursive orthogonal polynomials and fourth-order differential
  equation.
\newblock {\em J. Comput Appl. Math.}, 25 (1):$\ $105--109, 1989.

\bibitem{SHO50}
J.~A. Shohat and J.~D. Tamarkin.
\newblock The problem of moments.
\newblock In {\em Mathematical Surveys, volume 1}. American Mathematical
  Society, Providence, Rhode Island, 1950.

\bibitem{SLA66}
L.~J. Slater.
\newblock {\em Generalized Hypergeometric Functions}.
\newblock Cambridge University Press, Cambridge, 1966.

\bibitem{WAT44}
G.~N. Watson.
\newblock {\em Theory of {Bessel} functions}.
\newblock Cambridge University Press, Cambridge, 1944.

\bibitem{WIM87}
J.~Wimp.
\newblock Explicit formulas for the associated {Jacobi} polynomials and some
  applications.
\newblock {\em Canadian J. Math.}, 39:$\ $ 983--1000, 1987.

\end{thebibliography}

\enddocument